# EXACT SIMULATION OF DIFFUSIONS

By Alexandros Beskos[1] and Gareth O. Roberts

*Lancaster University*

We describe a new, surprisingly simple algorithm, that simulates *exact* sample paths of a class of stochastic differential equations. It involves rejection sampling and, when applicable, returns the location of the path at a random collection of time instances. The path can then be completed without further reference to the dynamics of the target process.

**1. Introduction.** Exact simulation of stochastic differential equations (SDEs) is a notorious problem within the applied probability community. The objective of this paper is to present a first step toward the solution of this problem. We describe a new algorithm we call the Exact Algorithm for simulating a class of SDEs. It involves rejection sampling and, when applicable, returns *exact* draws from any finite-dimensional distribution of the solution of the SDE.

Let $B = \{B_t; 0 \le t \le T\}$ be a scalar Brownian motion. Consider the general type of the one-dimensional Itô diffusion:

$$(1) \qquad dX_t = b(X_t)\,dt + \sigma(X_t)\,dB_t, \qquad 0 \le t \le T, X_0 = x \in \mathbf{R}$$

for drift coefficient $b: \mathbf{R} \mapsto \mathbf{R}$ and diffusion coefficient $\sigma: \mathbf{R} \mapsto \mathbf{R}$. Under certain regularity conditions on $b$ and $\sigma$ it can be shown that (1) has a solution $\{X_t; 0 \le t \le T\}$ weakly unique, that is, all the solutions have identical finite-dimensional distributions. Weak uniqueness, relatively more general than pathwise uniqueness, is sufficient for simulation purposes. For a formal definition of (1) see, for instance, [7].

Mathematical models of this kind are used to describe the evolution of stochastic phenomena in a wide range of disciplines and most times it is

Received February 2004; revised September 2004.
[1]Supported by the Greek State Scholarsips Foundation.

important to be able to simulate them. In rare cases (1) has an explicit solution with identifiable transition density. However, in most of the models used in practice it is necessary to resort to a numerical solution of (1). This has traditionally implied the use of some of the time discrete approximation methods (Euler, Taylor's expansion, etc.) which rely on small time approximate increment distributions for the diffusion (for a detailed account of these methods see [5]). In many cases throughout this paper the results of our algorithm will be compared with the Euler scheme which approximates (1) via the recursion

$$X_{t+\Delta} = X_t + b(X_t)\Delta + \sigma(X_t) \cdot \mathcal{N}(0, \Delta)$$

where we denote by $\mathcal{N}(\mu, \Sigma)$ the Gaussian distribution with mean $\mu$ and variance $\Sigma$.

For a large class of processes the Exact Algorithm provides an alternative which involves no approximation (apart from that inherent in any computer simulation) and yet is computationally highly efficient. It returns skeletons of *exact* paths which can be easily *filled in* without further reference to the diffusion dynamics. Thus the method can be used to simulate the diffusion at a prescribed collection of time points, or alternatively at times which occur to be interesting to the user, *after* the completion of the algorithm.

We begin (Section 2) by stating the rejection sampling technique in a way that serves our purposes. In Section 3 we present our method and in Section 4 give some results related with its efficiency. In Section 5 we apply the algorithm to a specific SDE and in Section 6 we take advantage of the properties of the Exact Algorithm to simulate exactly extremes and hitting times for the same SDE. Finally (Section 7), we present some ideas that could, in future research, overcome the restrictions of the algorithm and give some general conclusions. We are going to restrict the exposition of our algorithm to the case when $\sigma = 1$. This is not by any means restrictive since the SDE (1) can be transformed into one of unit diffusion coefficient for the process $Y = \{Y_t; 0 \le t \le T\}$ defined as

$$Y_t = \int_z^{X_t} \frac{1}{\sigma(u)}\, du,$$

where $z$ is an arbitrary element of the state space of $X$.

**2. A general rejection sampling algorithm.** Rejection sampling (RS) is a widely used simulation technique. It is frequently presented as follows. Assume that $f, g$ are probability densities w.r.t. some measure on $\mathbf{R}^d$ and that there exists $\varepsilon > 0$ such that $\varepsilon \frac{f}{g} \leq 1$. Then the iterative algorithm:

REJECTION SAMPLING
1. SAMPLE $Y \sim g$.
2. SAMPLE $U \sim \text{Unif}(0,1)$.
$3_1$. IF $U < \varepsilon \frac{f}{g}(Y)$ RETURN $Y$.
$3_2$. ELSE GO TO 1.

returns an observation distributed according to $f$.

This traditional rejection sampling algorithm seems to imply a fixed order in the acquisition of the required random elements: specifically, the proposed variate $Y$ precedes the decision variate $U$. However, such a prescribed order for the simulation steps is not necessary. The algorithm to be presented in Section 3 is more easily understood when imagining that the proposed variate actually *succeeds* the decision variate. Moreover, this reordering of the random steps can be carried out without adding any complexity to the algorithm. See [9] for a case where a similar reordering of the random inputs is used to perform an otherwise impossible MCMC algorithm for the Dirichlet mixture model.

Additionally, in some cases there are ways of constructing a condition for the acceptance or the rejection of the current proposed element $Y$ from minimal information about it. This will be essential in the diffusion context where it will never be possible to store a complete, continuous path. A similar idea appears in [4] on a perfect simulation algorithm of point processes. In that case a probability related with a complicated surface is expressed as the probability of an appropriately constructed event whose truth or otherwise is easy to verify.

We now present a formal definition of the RS algorithm in a way that incorporates the observations just described. Let $(S, \mathcal{S})$ be a sufficiently regular measurable space and $\nu, \mu$ probability measures on it such that $\mu$ is absolutely continuous w.r.t. $\nu$. Assume that there exists $\varepsilon > 0$ such that $f := \varepsilon \frac{d\mu}{d\nu} \leq 1$ $\nu$-a.s. and that it is easy to sample from $\nu$. The following proposition can be used to return draws from $\mu$.

PROPOSITION 1 (Rejection sampling). *Let $(Y_n, I_n)_{n \geq 1}$ be a sequence of i.i.d. random elements taking values in $S \times \{0, 1\}$ such that $Y_1 \sim \nu$ and $P[I_1 = 1 | Y_1 = y] = f(y)$ for all $y \in S$. Define $\tau = \min\{i \geq 1 : I_i = 1\}$. Then $P[Y_\tau \in dy] = \mu(dy)$.*

For the proof see the Appendix.

This presentation of the RS scheme does not assume any order for the simulation of $Y$ and $I$ and, besides the certain conditional property given in the proposition, does not restrict in any other way the construction of the binary indicator $I$.

**3. The Exact Algorithm.** Consider the stochastic process $X = \{X_t; 0 \leq t \leq T\}$ determined as the unique solution of the SDE

$$dX_t = \alpha(X_t)\, dt + dB_t, \qquad 0 \leq t \leq T, X_0 = 0. \tag{2}$$

The drift function $\alpha : \mathbf{R} \mapsto \mathbf{R}$ is presumed to satisfy the regularity properties that guarantee the existence of a global, weakly unique solution for (2). In particular, it suffices that $\alpha$ is locally Lipschitz; that is, for each $M > 0$ there exists $K_M > 0$ such that

$$|\alpha(y) - \alpha(x)| \leq K_M |y - x|; \qquad |x| \leq M, |y| \leq M$$

with a linear growth bound; there exists $K > 0$ such that

$$|\alpha(x)|^2 \leq K^2(1 + x^2), \qquad x \in \mathbf{R}.$$

See Chapter 4 of [5] for a detailed presentation of weaker conditions.

Since we are going to carry out rejection sampling it is convenient to think of the stochastic processes involved as measures induced on the space $C$ of continuous functions from $[0, T]$ to $\mathbf{R}$. We denote by $\omega$ a typical element of this space. Consider the coordinate functions $B_t(\omega) = \omega(t)$, $t \in [0, T]$, and the $\sigma$-field $\mathcal{C} = \sigma(\{B_t; 0 \leq t \leq T\})$. To avoid confusion between the coordinate functions and the Brownian motion process, we use the generic notation $BM = \{BM_t; 0 \leq t \leq T\}$ for a Brownian motion (BM) started at 0. Let $\mathbb{W}$ be the Wiener measure on $(C, \mathcal{C})$ so that $B = \{B_t; 0 \leq t \leq T\}$ is a BM.

3.1. *Rejection sampling for SDEs.* We explain how a rejection sampling algorithm can be set up for the case of the SDE given in (2). The analysis begins with the Girsanov transformation of measures. Let $\mathbb{Q}$ be the probability measure induced on $(C, \mathcal{C})$ by $X = \{X_t; 0 \leq t \leq T\}$.

PROPOSITION 2 (Girsanov transformation). *Assume that the drift coefficient $\alpha$ satisfies Novikov's condition:*

$$\mathrm{E}_{\mathbb{W}}\left[\exp\left\{\tfrac{1}{2}\int_0^T \alpha^2(B_t)\, dt\right\}\right] < \infty.$$

*It is then true that*

$$\frac{d\mathbb{Q}}{d\mathbb{W}}(\omega) = \exp\left\{\int_0^T \alpha(B_t)\, dB_t - \frac{1}{2}\int_0^T \alpha^2(B_t)\, dt\right\} =: G(B). \tag{3}$$

PROOF. See, for instance, [8], Chapter 8. □

Our goal is to implement a rejection sampling algorithm using (3) to construct an accept–reject mechanism. The difficulty is that exact evaluation of $G(B)$ is impossible. We can simplify $G(B)$ using Itô's formula to remove the Itô integral term and we shall see that this allows us to carry out the rejection scheme indirectly. We now need our first assumption.

CONDITION 1. *The drift coefficient $\alpha$ is everywhere differentiable.*

Under Condition 1, $G(B)$ admits the following simplification: let $A(u) = \int_0^u \alpha(y)\, dy$, $u \in \mathbf{R}$. Itô's formula then gives that

$$\int_0^T \alpha(B_t)\, dB_t = A(B_T) - A(B_0) - \tfrac{1}{2}\int_0^T \alpha'(B_t)\, dt.$$

We can now write $G(B)$ as

$$G(B) = \exp\Big\{A(B_T) - A(B_0) - \tfrac{1}{2}\int_0^T (\alpha^2(B_t) + \alpha'(B_t))\, dt\Big\}.$$

Rejection sampling using Brownian candidates is only conceivably possible if $G(B)$ is almost surely bounded and this is likely to require $A$ to be bounded. To remove this requirement we introduce a third probability measure which will be used to construct the candidates for the rejection sampling scheme.

Consider the *biased* Brownian motion $\overline{BM} = \{\overline{BM}_t; 0 \leq t \leq T\}$ heuristically defined as $(BM | BM_T = \rho)$ with $\rho$ distributed according to some density function $h : \mathbf{R} \mapsto [0, \infty)$ w.r.t. the Lebesgue measure.

PROPOSITION 3 (Biased Brownian motion). *Let $\mathbb{Z}$ be the probability measure induced by $\overline{BM}$ on $(C, \mathcal{C})$. If the support of $h$ is the real line, then $\mathbb{Z}$ is equivalent to $\mathbb{W}$ and*

$$\frac{d\mathbb{Z}}{d\mathbb{W}}(\omega) = \frac{h(B_T)}{(1/\sqrt{2\pi T})\exp(-B_T^2/(2T))}.$$

For the proof see the Appendix.

It is now trivial that

$$\frac{d\mathbb{Q}}{d\mathbb{Z}}(\omega) = \frac{d\mathbb{Q}}{d\mathbb{W}}\frac{d\mathbb{W}}{d\mathbb{Z}}(\omega)$$

$$\propto \exp\Big\{A(B_T) - \frac{B_T^2}{2T} - \frac{1}{2}\int_0^T (\alpha^2(B_t) + \alpha'(B_t))\, dt\Big\}\Big/ h(B_T),$$

where $\propto$ implies that we omitted some factors not depending on $\omega$.

CONDITION 2. $\int_{\mathbf{R}} \exp\{A(u) - u^2/2T\}\, du =: c < \infty.$

Under Condition 2 and after choosing $h(u) = \exp\{A(u) - u^2/2T\}/c$,

(4) $$\frac{d\mathbb{Q}}{d\mathbb{Z}}(\omega) \propto \exp\Big\{-\int_0^T \Big(\frac{1}{2}\alpha^2(B_t) + \frac{1}{2}\alpha'(B_t)\Big) dt\Big\}.$$

Assume now that the functional involved in the above integral is bounded:

CONDITION 3. There exist constants $k_1, k_2 \in \mathbf{R}$ such that $k_1 \leq \frac{1}{2}\alpha^2(u) + \frac{1}{2}\alpha'(u) \leq k_2$ for any $u \in \mathbf{R}$.

We can then write (4) as

$$\frac{d\mathbb{Q}}{d\mathbb{Z}}(\omega) \propto \exp\left\{-\int_0^T \phi(B_t)\,dt\right\} \tag{5}$$

for a function $\phi \geq 0$ defined as $\phi(u) = \frac{1}{2}\alpha^2(u) + \frac{1}{2}\alpha'(u) - k_1$, $u \in \mathbf{R}$. This creates the possibility of performing rejection sampling with candidates from $\mathbb{Z}$ in order to sample from $\mathbb{Q}$, which is the objective. We can choose the length $T > 0$ of the time interval under consideration so that

$$0 \leq \phi(u) \leq T^{-1} \qquad \text{for any } u \in \mathbf{R}. \tag{6}$$

Just consider any $T \leq 1/(k_2 - k_1) = 1/R$ where $R$ is the *identifiable* range of $(\alpha^2 + \alpha')/2$ in the sense that the lower (upper) bound we can actually obtain for $(\alpha^2 + \alpha')/2$ can be less (greater) than its maximum (minimum) value. From this point on it is assumed that $T$ has been fixed so that (6) is true.

We remark that Condition 3 is trivially implied by

CONDITION 3'. There exist constants $k_1', k_2' \in \mathbf{R}$ such that $k_1' \leq \alpha(u)$, $\alpha'(u) \leq k_2'$ for any $u \in \mathbf{R}$.

Furthermore, Condition 3' implies that Condition 2 holds and moreover the constants $k_1'$ and $k_2'$ point the way to implementing simple rejection sampling algorithms for simulating the endpoint of the biased Brownian motion.

3.2. *Constructing the Exact Algorithm.* We set $H(\omega) = \int_0^T \phi(B_t)\,dt$. We assume from now on that we have access to paths $\omega \sim \mathbb{Z}$ of the biased Brownian motion $\overline{BM}$. Note that we can realize such a path at any finite collection of time instances by drawing first its ending point $\omega_T \sim h$ and then the rest of the skeleton according to the dynamics of a Brownian bridge. Drawing from the univariate distribution $h$ cannot be a big problem; [1] gives many algorithms for drawing from densities on $\mathbf{R}$ that could be used as envelopes for a rejection sampling scheme on $h$.

The preliminaries of Section 3.1 together with the general rejection sampling protocol of Proposition 1 ensure that the following algorithm (were it implementable in practice) would output realizations of the diffusion $X$ that

solves the SDE (2):

IMPOSSIBLE ALGORITHM
1. SAMPLE A COMPLETE, CONTINUOUS PATH OF $\overline{BM}$, $\omega \sim \mathbb{Z}$.
2. COMPUTE $H(\omega)$.
3. PRODUCE A BINARY INDICATOR $I$ SUCH THAT $\mathrm{P}[I=1|\omega] = \exp\{-H(\omega)\}$.
4. IF $I=0$ GO TO 1.
5. OUTPUT $\omega$.

The indicator $I$ is easily constructed with the use of some $U \sim \mathrm{Unif}(0,1)$. In practice, we can only simulate the path $\omega$ at any given finite collection of instances $0 \leq t_1, t_2, \ldots, t_n \leq T$ so evaluating the integral $H(\omega)$ is impossible. However, given that Conditions 1–3 are satisfied, we can produce an algorithm which manages to circumvent steps 1 and 2 and still carry out steps 3 and 4 exactly given only a *finite* but random skeleton of instances of $\omega$.

The idea builds on the simple observation that for a bounded function $0 \leq \phi(u) \leq T^{-1}$ events of probability $\int_0^T \phi(u)\,du$ can be constructed simply by drawing a random point $(V,W) \sim \mathrm{Unif}[(0,T) \times (0,T^{-1})]$. Then the event $\{\phi(V) \geq W\}$ will have the required probability.

To extend this idea to an event of probability $\exp(-H)$ we exploit a Taylor series expansion construction which gives us an event of probability $\exp(-H)$ as an event which depends upon a countable sequence of events of probability $H$. Furthermore, it turns out to be possible to express this event both as the countable union of a sequence of increasing events and as the countable intersection of another sequence of decreasing events. Crucially, for each event in either of the two sequences its truth or otherwise can be confirmed by a finite skeleton of a path $\omega \sim \mathbb{Z}$ and, consequently, the truth or otherwise of the event of probability $\exp\{-H(\omega)\}$ can also be determined after finite computations.

All the above ideas are presented in a rigorous way in Theorem 1 that follows. The construction to be described is similar in spirit, though in a different context, with Von Neumann's comparison method for the simulation of exponential random variables; see [2] for a detailed review. We have denoted by $(\Omega, \mathcal{F}, \mathrm{Prob})$ the underlying probability space that generates all the random elements involved in the theorem.

THEOREM 1. *Let $\omega \sim \mathbb{Z}$ be a path of the biased Brownian motion $\overline{BM}$ on $[0,T]$. Let $\tau = (V_i, W_i)_{i \geq 1}$ be a sequence of i.i.d. points uniformly distributed on $(0,T) \times (0, 1/T)$. Consider also some $U \sim \mathrm{Unif}(0,1)$. Assume that $\omega$, $\tau$ and $U$ are independent. Define the following events:*

$$\Gamma_0 = \Omega, \qquad \Gamma_n = \left\{ \phi(B_{V_1}(\omega)) \geq W_1, \ldots, \phi(B_{V_n}(\omega)) \geq W_n, U \leq \frac{1}{n!} \right\},$$
(7)
$$n = 1, 2, \ldots.$$

*Consider the sequence of events* $(E_n)_{n\geq 1}$ *defined as*

$$E_{2n+1} = (\Gamma_0 - \Gamma_1) + (\Gamma_2 - \Gamma_3) + \cdots + (\Gamma_{2n} - \Gamma_{2n+1}),$$
$$n = 0, 1, \ldots,$$

(8)
$$E_{2n+2} = (\Gamma_0 - \Gamma_1) + (\Gamma_2 - \Gamma_3) + \cdots + (\Gamma_{2n} - \Gamma_{2n+1}) + \Gamma_{2n+2},$$
$$n = 0, 1, \ldots,$$

*where* $(+)$ *implies union of disjoint sets and* $D - F = D \cap F^c$ *for any sets* $F \subseteq D$. *Then:*

(i) $(E_{2n+1})_{n\geq 0}$, $(E_{2n+2})_{n\geq 0}$ *are sequences of increasing and decreasing events, respectively, with* $E_{2\kappa+1} \subseteq E_{2\lambda+2}$ *for any* $\kappa, \lambda \in \{0, 1, \ldots\}$ *and* $\text{Prob}[\bigcap_0^\infty E_{2n+2} - \bigcup_0^\infty E_{2n+1}] = 0$.

(ii) *Let* $E = \bigcup_0^\infty E_{2n+1}$. *If* $I$ *is a binary indicator such that* $I = 1$ *when* $E$ *occurs and* $0$ *otherwise, then*

$$\text{Prob}[I = 1|\omega] = \exp\left\{-\int_0^T \phi(B_t(\omega))\,dt\right\}.$$

See the proof in the Appendix.

Recall that $B_t$ is the coordinate mapping $B_t(\omega) = \omega(t)$, $t \in [0, T]$. We will from now on write $\omega(t)$ instead of $B_t(\omega)$. In practice, we can simply identify $\bigcap_0^\infty E_{2n+2} - \bigcup_0^\infty E_{2n+1}$ with the null event $\varnothing$ so that $E = \bigcup_0^\infty E_{2n+1} = \bigcap_0^\infty E_{2n+2}$. Figure 1 illustrates the decomposition of $\Omega$ over the sets $E_i$, $i \geq 1$.

Theorem 1 does not impose any restriction on the order of the realization of the random elements when rejection sampling will take place. To transform the theorem into a feasible rejection sampling algorithm it is necessary that we draw $U$ and then generate the path $\omega \sim \mathbb{Z}$ and the sequence $\tau$ in parallel. We follow an iterative process that involves drawing $(V_i, W_i)$ uniformly from $(0, T) \times (0, 1/T)$ and then simulating $\omega(V_i)$ conditionally on the already obtained $\omega(V_1), \omega(V_2), \ldots, \omega(V_{i-1})$, for $i \geq 1$. Recall that this recursive construction of the path of $\overline{BM}$ is straightforward as long as someone begins by simulating its ending point $\omega(T) \sim h$ (the choice for $h$ is given after Condition 2 of Section 3.1). Then, the rest of the path is a Brownian bridge; given that the locations $\{\omega(V_1), \ldots, \omega(V_{i-1}), \omega(T)\}$ have been constructed the path can be realized at the instance $V_i$ just by drawing

$$\omega(V_i) \sim \mathcal{N}\left(\omega(V_-) + \frac{\omega(V_+) - \omega(V_-)}{V_+ - V_-}(V_i - V_-), \frac{(V_+ - V_i)(V_i - V_-)}{V_+ - V_-}\right),$$

where we have defined $V_- = \max\{0, V_j, j = 1, \ldots, i-1 : V_j < V_i\}$ and $V_- = \min\{T, V_j, j = 1, \ldots, i-1 : V_j > V_i\}$. See, for instance, page 360 of [3] for the representation of a Brownian bridge as a transformation of an unconditional Brownian motion which implies the above formula.

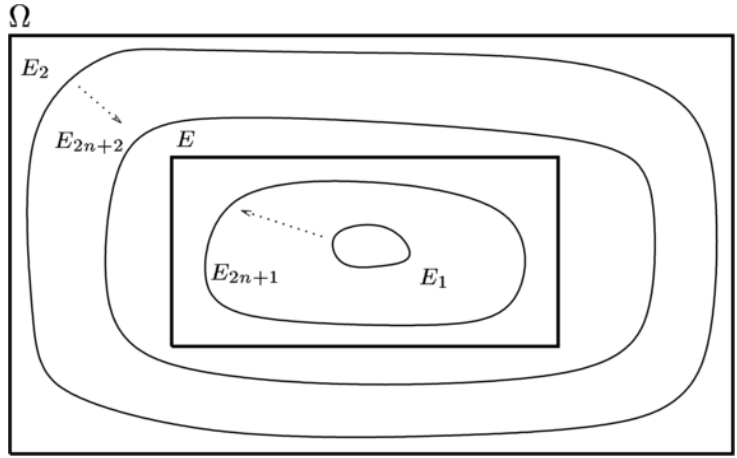

FIG. 1. *An illustration of the sets* $E_i, i \geq 1$. *When some even-numbered* $E_i$ *does not occur or an odd-numbered* $E_i$ *occurs,* $E$ *does not occur or occurs, respectively.*

From the definition of $(E_i)_{i \geq 1}$ it is clear that after $j$ iterations are carried out we have the necessary information to decide if any of the events $E_1, E_2, \ldots, E_j$ occurred or not. We perform iterations until the first time that an odd-numbered $E_i$ occurs or an even-numbered $E_i$ does not occur. It is then clear from Theorem 1 that the former case gives $I = 1$ and the latter $I = 0$. Further realization of the random elements will not change the decision about $I$. If needed, we can continue simulating an accepted path of $\overline{BM}$ so that the resulted path of $X$ is constructed at any time instances requested.

The recursive definition of $(E_i)_{i \geq 1}$ allows for a simple way of carrying out the successive steps of the iterative method described above. Assume for instance that $2n$ iterations have taken place without a decision about $I$, that is, $E_{2n}$ happened and $E_{2n-1}$ did not happen. From (7) and (8) it is clear that

$$E_{2n} = E_{2n-1} + \Gamma_{2n} \quad \text{and}$$
$$E_{2n+1} = E_{2n} - \Gamma_{2n+1} = E_{2n-1} + (\Gamma_{2n} - \Gamma_{2n+1}).$$

The first equation indicates that $\Gamma_{2n}$ occurred and the second that we only need to check if the subevent of $\Gamma_{2n}$

$$\Gamma_{2n} - \Gamma_{2n+1} = \Gamma_{2n} \cap \left\{\phi(\omega(V_{2n+1})) < W_{2n+1} \text{ or } U > \frac{1}{(2n+1)!}\right\}$$

occurred [i.e., $\phi(\omega(V_{2n+1})) < W_{2n+1}$ or $U > \frac{1}{(2n+1)!}$] or not to reach to a similar conclusion about $E_{2n+1}$. The same convenient interpretation appears for the case when an even-numbered iteration, say the $2n$th, is carried out. Given that $I$ is not determined before that step we can conclude that $E_{2n}$ did not take place if $U > 1/(2n)!$ or $\phi(\omega(V_{2n})) < W_{2n}$.

We can now present the pseudocode that implements Theorem 1 to carry out rejection sampling:

EXACT ALGORITHM
1. INITIATE A PATH OF $\overline{BM}$: SET $\omega(0) = 0$ AND DRAW $\omega(T) \sim h$.
2. DRAW $U \sim \text{Unif}(0, 1)$. SET $i = 0$.
3. DRAW $(V, W) \sim \text{Unif}[(0, T) \times (0, 1/T)]$. SET $i = i + 1$.
4. CONSTRUCT $\omega(V)$ GIVEN THE CURRENTLY SIMULATED INSTANCES OF $\omega$.
$5_1$. IF $\phi(\omega(V)) < W$ OR $U > 1/i!$ THEN
   IF $i$ IS EVEN SET $I = 0$ AND GO TO 1.
   IF $i$ IS ODD SET $I = 1$ AND GO TO 6.
$5_2$. ELSE GO TO 3.
6. OUTPUT THE CURRENTLY SIMULATED INSTANCES OF $\omega$.

As already mentioned, we can interpose a step that constructs the proposed path at any time instances we require. It is now clear that the algorithm returns exact skeletons of the target process $X$ at any given finite collection of time instances.

3.3. *A factorization for* $\mathbb{Q}$. From a probabilistic point of view, the Exact Algorithm manages to decompose the target measure $\mathbb{Q}$ in terms of a product of Brownian bridges after appropriately extending the underlying probability space.

Analytically, let $\mathsf{S}$ be the random skeleton produced by the Exact Algorithm [we do not include in $\mathsf{S}$ the starting point $(0, 0)$]. If $S_k$ is the space of the possible configurations of $k$ points $\{x_1, x_2, \ldots, x_k\} \subseteq [0, T] \times \mathbf{R}$, $k \geq 1$, then $\mathsf{S}$ takes values in $\bigcup_k S_k$. We denote by $\mathsf{s} = \{(u_1, y_1), (u_2, y_2), \ldots, (u_k, y_k)\}$ a typical element of this state space. We avoid details about the $\sigma$-algebra construction on $\bigcup_n S_n$, and simply denote by $\mathcal{L}_\mathsf{S}$ the distribution of $\mathsf{S}$. Let $\mathbb{BB}(s, x; t, y)$, for $0 \leq s < t$, $x, y \in \mathbf{R}$, be the probability measure corresponding to a Brownian bridge starting at the time instance $s$ from $x$ and finishing at the time instance $t$ at $y$. In terms of probability measures, the Exact Algorithm manages to factorize $\mathbb{Q}$ in the following way:

$$(9) \qquad \mathbb{Q} = \bigotimes_{i=1}^{k} \mathbb{BB}(u_{i-1}, y_{i-1}; u_i, y_i) \otimes \mathcal{L}_\mathsf{S}(d\mathsf{s}),$$

where $(u_0, y_0) \equiv (0, 0)$. Critically, the rejection sampling construction of the Exact Algorithm allows for the simulation of $\mathcal{L}_\mathsf{S}$, so drawing from $\mathbb{Q}$ is then straightforward. Figure 2 gives a graphical illustration of the factorization (9). Because of this decomposition of $\mathbb{Q}$ it is possible to identify characteristics of the target process $X$ after considering the properties of the Brownian bridges that fill in its skeletons. Thus, it is possible to simulate exactly hitting times, extremes (we present these applications analytically in Section

6) and any other random elements for which there are explicit results for the case of Brownian paths.

Another simple example when we can exploit (9) is at the Monte Carlo evaluation of the expected value of functionals of $X_t$ for some time instance $t \in [0, T]$. Assume that $\{X_t^i, \mathsf{S}^i\}_{1 \leq i \leq n}$ are $n$ draws for the skeleton $\mathsf{S}$ and $X_t$ produced by the Exact Algorithm. An estimator of $\mathrm{E}[f(X_t)]$, for some function $f$, can be the mean of the $f(X_t^i)$, $1 \leq i \leq n$. Using the simple conditional expectation property $\mathrm{Var}[f(X_t)] \geq \mathrm{Var}[[\mathrm{E}[f(X_t)|\mathsf{S}]]$ we can produce an unbiased estimator of smaller variance after considering the mean of the $\mathrm{E}[f(X_t)|\mathsf{S}^i]$, $1 \leq i \leq n$. Conditionally on the skeleton, $X_t$ follows a normal distribution, so for reasonable functions $f$ finding the $\mathrm{E}[f(X_t)|\mathsf{S}]$ will be straightforward.

**4. Efficiency of the Exact Algorithm.** Two aspects of the Exact Algorithm need to be examined. The first involves the probability of accepting a proposed path of $\overline{BM}$ as a path of $X$. The second has to do with the number of points in $(0, T) \times (0, 1/T)$ required to reach a decision about a proposed path.

We can rewrite (5) as

$$\varepsilon(T) \frac{d\mathbb{Q}}{d\mathbb{Z}}(\omega) = \exp\left\{-\int_0^T \phi(B_t)\, dt\right\}$$

for some appropriate $\varepsilon(T)$. Note that $\mathbb{Q}$ and $\mathbb{Z}$ are both probability measures so $\varepsilon(T)$ equals precisely the probability that a path $\omega \sim \mathbb{Z}$ is accepted as a path from $\mathbb{Q}$ (see the proof of Proposition 1). Equivalently, $\varepsilon(T) = \mathrm{Prob}[I = 1]$ for the binary indicator $I$ defined in Theorem 1.

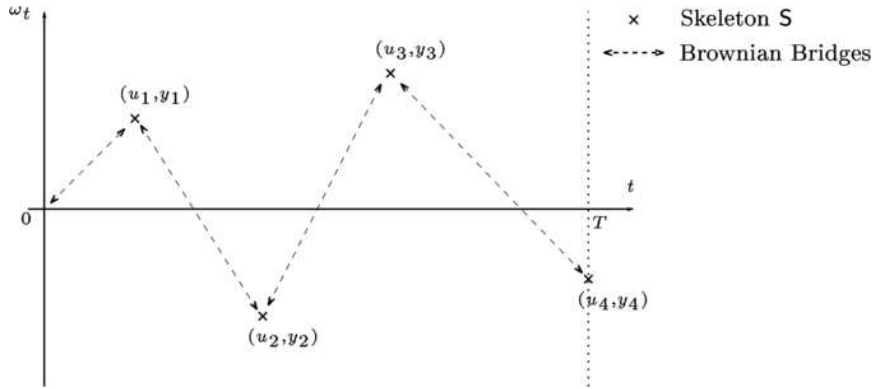

FIG. 2. *The factorization of $\mathbb{Q}$: drawing from $\mathbb{Q}$ is achieved after simulating the skeleton $\mathsf{S} = \{(u_1, y_1), (u_2, y_2), \ldots, (u_k, y_k)\}$ (in the case of the figure, $k = 4$) and then filling in the rest of the path with independent Brownian bridges.*

PROPOSITION 4. *For any appropriate $T$, the probability $\varepsilon(T)$ of the event $\{I=1\}$ for $I$ defined in Theorem 1 is at least $e^{-1}$ and is decreasing in $T$ with $\lim_{T\downarrow 0}\varepsilon(T)=1$.*

For the proof see the Appendix.

As expected, the probability of accepting the proposed paths of $\overline{BM}$ increases when we decrease the length of the time period under consideration.

It is important that we give a result for the number of points uniformly drawn in $(0,T)\times(0,1/T)$ required to accept or reject a proposed path. Recall that $(\Omega,\mathcal{F},\mathrm{Prob})$ is the underlying probability space that generates all the random elements involved in Theorem 1.

PROPOSITION 5. *Define the random variable $N:\Omega\mapsto\{1,2,\ldots\}$ as*

$$N(w)=\begin{cases}\min\{n=2,4,\ldots:w\notin E_n\}, & \text{if } w\in E^c,\\ \min\{n=1,3,\ldots:w\in E_n\}, & \text{if } w\in E.\end{cases}$$

*Then $\mathrm{E}[N]\leq e$.*

PROOF. Just note that $\mathrm{Prob}[N\geq n]\leq \mathrm{Prob}[U\leq \frac{1}{(n-1)!}]$, so it is straightforward that

$$\mathrm{E}[N]\leq\sum_{i=1}^{\infty}\frac{1}{(n-1)!}=e. \qquad\square$$

The random variable $N$ counts the number of events $E_n$ we have to consider before deciding if $w$ belongs to $E$ or $E^c$ (equivalently, if the corresponding realizations of the random elements involved in Theorem 1 make $E$ happen or not). It is a perhaps surprising result that for any eligible $T$ and any drift coefficient $\alpha$ we are expecting on average less than three of the points drawn uniformly from $(0,T)\times(0,1/T)$ before we decide for the acceptance or the rejection of a proposed path.

Note that the Exact Algorithm can take advantage of the Markov property of the process $X$ and produce skeletons of any requested length $l>0$ after merging skeletons of lengths acceptable by the algorithm. The choice of the length of the merged skeletons and, subsequently, the efficiency of the Exact Algorithm depend on the identifiable range of the functional $(\alpha^2+\alpha')/2$ of the drift function. Recall that we can identify analytically which satisfies

$$\sup_{x\in\mathbf{R}}(\alpha^2+\alpha')(x)/2-\inf_{x\in\mathbf{R}}(\alpha^2+\alpha')(x)/2\leq R.$$

The following proposition gives a result for the case when the Exact Algorithm merges skeletons of the maximum eligible length $T=1/R$ to obtain a skeleton of length $l$. $\lceil u\rceil$ is the minimum integer not smaller than $u\in\mathbf{R}$.

PROPOSITION 6. *If $R$ is the identifiable range of the functional $(\alpha^2 + \alpha')/2$ of the drift $\alpha$ and $N_l$ is the total number of the uniformly drawn points needed for the Exact Algorithm to return a skeleton of length $l > 0$, then*

$$\mathrm{E}[N_l] \leq \lceil l \cdot R \rceil \times e^2.$$

For the proof see the Appendix.

In total, the Exact Algorithm requires elementary computational skills. Except for the appealing characteristic of being exact it seems to compete with conventional approximation techniques even in terms of time efficiency. The example that follows favors this assertion.

**5. Applying the Exact Algorithm.** We apply the Exact Algorithm to the SDE:

(10) $$dX_t = \sin(X_t)\, dt + dB_t, \qquad 0 \leq t \leq T, X_0 = 0.$$

Conventional methods can only approximate sample paths for the solution $X$ of (10) after resorting to one of the suggested time discretization techniques; the Exact Algorithm returns exact skeletons of $X$.

The drift coefficient $\alpha \equiv \sin$ satisfies Conditions 1 and 2. It is also easy to check that $-1/2 \leq \frac{1}{2}\sin^2(u) + \frac{1}{2}\cos(u) \leq 5/8$ for any $u \in \mathbf{R}$, so Condition 3 is satisfied for $k_1 = -1/2$ and $k_2 = 5/8$. In the present context $\phi(u) = \frac{1}{2}\sin^2(u) + \frac{1}{2}\cos(u) + \frac{1}{2}$. We can now choose $T = 1/(k_2 - k_1) = 8/9$; for this ending time instance, $0 \leq \phi(u) \leq T^{-1}$ for any real $u$. The $\overline{BM}$ process is defined w.r.t. a BM via the rule $\overline{BM} = (BM | BM_T = \rho)$ with $\rho \sim h \propto \exp\{-\cos(u) - u^2/2T\}$. We can draw efficiently from this univariate distribution using rejection sampling with Gaussian proposals.

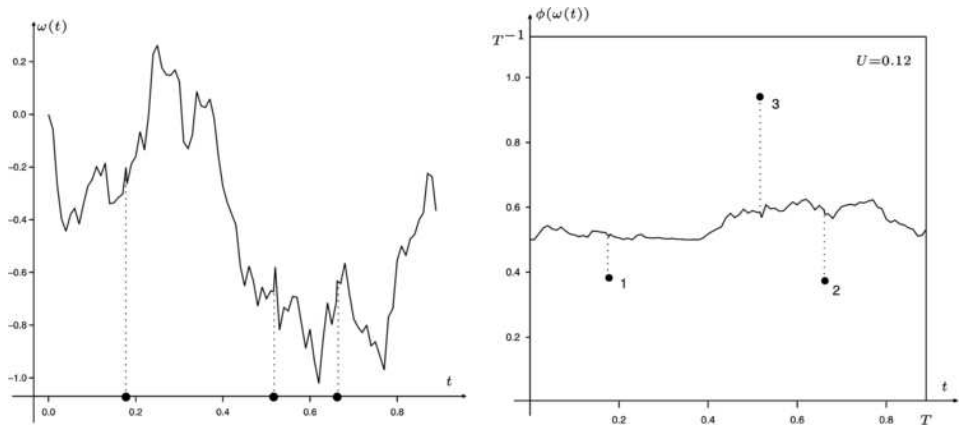

FIG. 3. *A case when the Exact Algorithm accepts the proposed path $\omega \sim \mathbb{Z}$.*

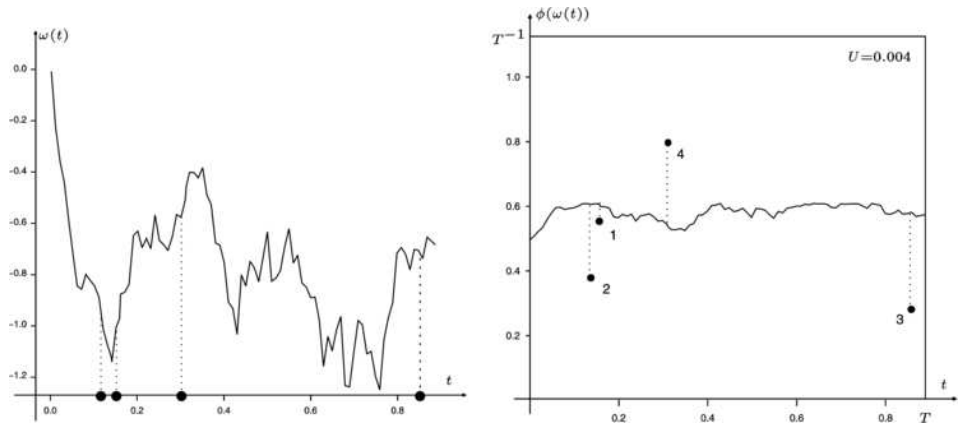

Fig. 4. *A case when the Exact Algorithm rejects a proposed path $\omega \sim \mathbb{Z}$.*

Everything is now set up for applying the Exact Algorithm. We run the algorithm until we generate 5000 exact skeletons of $X$. We had to propose 12,320 paths of $\overline{BM}$ to get the exact paths so we can estimate that Prob$[I=1]$ for the indicator of Theorem 1 is close to 0.41. In 58% of the proposed paths the decision about accepting or rejecting a path was taken after simulating 1 and 2 points respectively uniformly from $(0,T) \times (0, T^{-1})$. The maximum number of points needed for the acceptance or the rejection of a path was 7 and 6, respectively.

Figure 3 shows on the left an exact skeleton of $X$ and on the right the same skeleton after considering the transformation $x \mapsto \phi(x)$ for each of its joints. In the case of the $\phi$-path the square encloses the area $(0,T) \times (0, T^{-1})$ and the black spots show the location of three points uniformly drawn from this rectangle. The numbers next to each circle show the order with which they were obtained. At the top right corner we have written the draw $U \sim$ Unif$(0,1)$ needed by the Exact Algorithm. The third point exceeded the graph of $\phi$ so the algorithm decided that the event $E_3$ of Theorem 1 occurred after realizing the proposed path $\omega$ only at three time instances. The rest of the path of $X$ involves the time instances $0.01i$ for all $i=1,2,\ldots$ such that $0.01i < T$ and is produced after filling in the proposed path.

Figure 4 shows similar graphs for the case when a proposed path is rejected. The graph that involves $\phi$ indicates that the proposed path was rejected because the event $E_4$ did not occur.

Figure 5 shows the estimated density for the distribution of $X_1$ from samples of size 1,000,000 obtained after using the Exact Algorithm and the Euler approximation method for different time discretizations. For the case of the Exact Algorithm, to draw from $X_1$ we had to merge skeletons on the time intervals $[0, 8/9]$ and $[8/9, 1]$. Table 1 presents the times in seconds needed to get these samples and the $p$-values of the Kolmogorov–Smirnov test that compares the approximate draws of the Euler method with the exact ones of our algorithm. All programs were written in C-language and executed on a Athlon PC, running at 1500 MHz. It is clear, at least for the case of the SDE given in (10), that the Exact Algorithm is superior to the Euler approximation even in terms of computational time.

Note also that, as implied in Proposition 6, the time needed to draw from $X_l$, $l > 0$, increases linearly in $l$ for the Exact Algorithm. In contrast, the Euler method needs thinner approximation to produce reliable results as $l$ increases because of the accumulating errors.

**6. Exact simulation of extremes and hitting times.** We take advantage of the properties of the Exact Algorithm to simulate exactly the maximum and the hitting time of a horizontal line boundary for one-dimensional diffusions determined by SDEs of the type (2) under the Conditions 1–3 given in Section 3.1.

As explained in Section 3.3, the Exact Algorithm reduces the problem of detecting characteristics for a path of the process $X$ to the much more straightforward task of carrying out the certain detection process for its Brownian subpaths. To simplify the presentation of the algorithms that follow we denote by $\mathcal{S}_l\{x; t_1, t_2, \ldots, t_n\}$ a skeleton of the target process $X$ at the time instances $0 < t_1 < t_2 < \cdots < t_n = l$ starting at $x$.

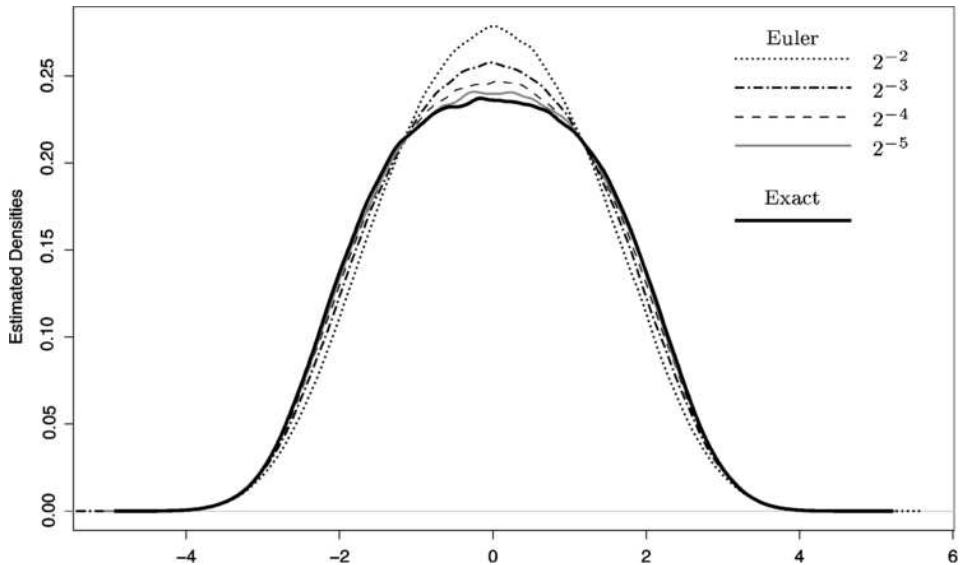

FIG. 5. *The estimated density for the distribution $X_1$ from samples of size* 1,000,000 *generated by the Euler approximation (four cases for time increments $2^{-2}$, $2^{-3}$, $2^{-4}$, $2^{-5}$) and the Exact Algorithm.*

TABLE 1
*Times in seconds needed for the Exact Algorithm and the Euler approximation to produce* 1,000,000 *draws from the distribution of* $X_1$

| Exact | Euler (for inc. $\Delta = 2^{-n}, n = 2, 3, \ldots, 8$) | KS test |
|---|---|---|
|  | 5 sec | 0 |
|  | 11 sec | 0 |
|  | 15 sec | 0 |
| 35 sec | 24 sec | 0 |
|  | 45 sec | 0.002 |
|  | 87 sec | 0.261 |
|  | 174 sec | 0.412 |

The last column shows the *p*-values for the Kolmogorov–Smirnov test with null hypothesis that the exact and the corresponding approximate draws come from the same distribution.

For the case of the maximum value of $X$ over the time interval $[0, l]$, $l > 0$, we apply the Exact Algorithm, if necessary after dividing $[0, l]$ in smaller pieces of length at most the maximum length $T$ permitted by the Exact Algorithm and dealing with each piece separately, and locate the maximum value of the accepted proposed paths of $\overline{BM}$ after drawing the maximum of all the Brownian bridges (BBs) that intervene between the successive unveiled instances of these paths. Drawing the maximum of a BB is straightforward.

Let $BM^y = \{BM^y_s; 0 \leq s \leq t\}$ be a Brownian motion over the interval $[0, t]$ started at $y$. If $M_t = \sup\{BM^y_s; s \in [0, t]\}$, then it is a known result that

$$P[M_t \in db | BM^y_t = a] \propto (2b - y - a) \exp\left\{-\frac{(2b - y - a)^2}{2t}\right\} db,$$
$$b \geq \max\{y, a\}.$$

For a proof, see, for instance, [3], page 95. This is just a linear transformation of a Rayleigh distribution and it is easy to verify that

$$[M_t | BM^y_t = a] \stackrel{d}{=} \tfrac{1}{2}(\sqrt{2tE(1) + (a-y)^2} + a + y),$$

where $E(1)$ denotes an exponential random variable with unit mean. This formula generates the maximum of a BB of length $t$ starting at $y$ and finishing at $a$ (for arbitrary $t$, $y$ and $a$) and will be applied for all the BBs that fill in the exact skeleton of $X$. Thus, the complete algorithm for drawing the maximum of $X = \{X_t; 0 \leq t \leq l\}$ is as follows:

EXACT SIMULATION OF THE MAXIMUM
1. CALL ON THE EXACT ALGORITHM AND DRAW A SKELETON
   $\mathcal{S}_l\{0; t_1, t_2, \ldots, t_n\}$.
2. SIMULATE THE MAXIMA $M^1, M^2, \ldots, M^n$ OF THE BBS THAT FILL IN $\mathcal{S}_l$.
3. OUTPUT $\sup\{M^i; i = 1, 2, \ldots, n\}$.

A similar procedure is carried out for the simulation of first passage times. We now have to check if each of the BBs that fill in the skeleton of $X$ produced by the Exact Algorithm hit some arbitrary boundary $\gamma$ or not and for the first bridge that hits $\gamma$ to find the precise time instance when that occurs. Assume for simplicity that $\gamma > 0$ is greater than the starting point $0$ of $X$.

The hitting time of a horizontal boundary for a BB is closely related with the first passage of an unconditional Brownian motion over a linear boundary. It is a known result that if $BM(\delta) = \{BM_s(\delta); 0 \leq s \leq t\}$ is a BB of length $t$ started at $0$ and finishing at $\delta$ and $BM = \{BM_s; s \geq 0\}$ is an unconditional Brownian motion started at $0$, then

$$(11) \qquad BM_s(\delta) \stackrel{d}{=} \frac{s}{t}\delta + \frac{t-s}{\sqrt{t}} BM_{s/(t-s)}, \qquad 0 \leq s \leq t.$$

Let $\tau_\gamma(\delta) = \inf\{s \in [0,t] : BM_s(\delta) \geq \gamma\}$ and $\tau_{\eta,\zeta} = \inf\{s \geq 0 : BM_s \geq \eta + \zeta s\}$ with the convention that $\inf \varnothing = \infty$. Then for any $s \in [0,t]$ it is true that

$$\begin{aligned}
\mathrm{P}[\tau_\gamma(\delta) > s] &= \mathrm{P}[BM_u(\delta) < \gamma, \text{ for all } 0 \leq u \leq s] \\
&= \mathrm{P}\left[\frac{u}{t}\delta + \frac{t-u}{\sqrt{t}} BM_{u/(t-u)} < \gamma, \text{ for all } 0 \leq u \leq s\right] \\
&= \mathrm{P}\left[BM_{u^*} < \frac{\gamma}{\sqrt{t}} + \frac{\gamma - \delta}{\sqrt{t}} u^*, \text{ for all } 0 \leq u^* \leq \frac{s}{t-s}\right] \\
&= \mathrm{P}\left[\tau_{\eta,\zeta} > \frac{s}{t-s}\right]
\end{aligned}$$

for $\eta = \gamma/\sqrt{t}$, $\zeta = (\gamma - \delta)/\sqrt{t}$. It is now clear that for these values of $\eta$, $\zeta$:

$$(12) \qquad \tau_\gamma(\delta) \stackrel{d}{=} g(\tau_{\eta,\zeta}) \qquad \text{for } g(u) = \begin{cases} tu/(u+1), & 0 < u < \infty, \\ \infty, & u = \infty. \end{cases}$$

The density of $\tau_{\eta,\zeta}$ is given by the Bachelier–Lévy formula:

$$(13) \qquad p^{\eta,\zeta}(u) = \frac{|\eta|}{u^{3/2}\sqrt{2\pi}} \exp\{-(\zeta u + \eta)^2/2u\}, \qquad u > 0.$$

A proof can be found in [6], Chapter 1. If we denote by $IG(\mu, \lambda)$, $\mu > 0$, $\lambda > 0$, the inverse Gaussian distribution with density:

$$IG(\mu, \lambda, u) = \sqrt{\frac{\lambda}{2\pi u^3}} \exp\left\{-\frac{\lambda(u-\mu)^2}{2\mu^2 u}\right\}, \qquad u > 0,$$

then, when $\eta\zeta < 0$ the density (13) is easily identified as an $IG(-\eta/\zeta, \eta^2)$ and integrates to 1, that is, the Brownian path hits the linear boundary $\eta + \zeta t$ with probability 1, as expected. When $\eta\zeta > 0$, (13) can be written as $\exp(-2\zeta\eta)IG(\eta/\zeta, \eta^2, u)$, which means that the Brownian path hits the linear boundary w.p. $\exp(-2\zeta\eta)$ and when it does the distribution of the hitting time is just the $IG(\eta/\zeta, \eta^2)$. Draws from the inverse Gaussian distribution can be generated in a very efficient way; see the algorithm described in [1], Chapter IV. Thus, simulating $\tau_{\eta,\zeta}$ and, via (12), $\tau_\gamma(\delta)$ for any values of the involved parameters is straightforward.

We now have all that is needed to carry out the algorithm for the exact simulation of the hitting time of a horizontal boundary for the solution $X$ of (2). Since it is not possible to know a priori the length of the path of $X$ we need to construct before we locate the time instance when $X$ hits the boundary, we merge as many exact skeletons of $X$ of eligible length $T$ as necessary before some intervening BB hits the boundary.

EXACT SIMULATION OF HITTING TIME
1. SET $i=1$, $y=0$.
2. CALL ON THE EXACT ALGORITHM AND GENERATE
   $\mathcal{S}_T\{y; t_{i,1}, t_{i,2}, \ldots, t_{i,n_i}\}$.
3$_1$. IF THE INTERVENING BBS HIT $\gamma$, SAVE THE FIRST TIME $\tau$ WHEN THAT OCCURS.
3$_2$. ELSE, SET $y = X_T$, $i = i+1$ AND GO TO 2.
4. RETURN $(i-1)T + \tau$.

Assume that $N^\gamma$ is the total number of the uniformly drawn points needed for finding $\tau_\gamma = \inf\{t \geq 0 : X_t = \gamma\}$ and that $N_i$ is the number of these points needed for accepting a proposed path from $(i-1)T$ to $iT$, $i \geq 1$. Let $\mathcal{G}_i$ be the information for the obtained exact path of $X$ until the time $iT$, $i \geq 0$. Proposition 6 is true for any given starting point of the target process, so $\mathrm{E}[N_i|\mathcal{G}_{i-1}] \leq e^2 \lceil T \cdot R \rceil$ for all $i \geq 1$, for $R$ as defined in Proposition 6. We set $\tau := \lceil \tau_\gamma/T \rceil$; $\tau$ is the total number of the exact skeletons we have to merge before we find $\tau_\gamma$. Clearly,

$$\mathrm{E}[N^\gamma] = \sum_{i=1}^\infty \mathrm{E}[N_i \cdot \mathbb{I}\{\tau \geq i\}] = \sum_{i=1}^\infty \mathrm{E}[\mathrm{E}[N_i \cdot \mathbb{I}\{\tau \geq i\}|\mathcal{G}_{i-1}]]$$
$$= \sum_{i=1}^\infty \mathrm{E}[\mathbb{I}\{\tau \geq i\} \cdot \mathrm{E}[N_i|\mathcal{G}_{i-1}]] \leq e^2\lceil T \cdot R\rceil \sum_{i=1}^\infty \mathrm{P}[\tau \geq i]$$
$$\leq e^2\lceil T \cdot R\rceil \left(\frac{\mathrm{E}[\tau_\gamma]}{T} + 1\right),$$

where we have used the fact that $\{\tau \geq i\}$ is $\mathcal{G}_{i-1}$-measurable. So the expected time for the termination of the algorithm, in terms of the uniformly drawn points needed, is finite when $\mathrm{E}[\tau_\gamma]$ is finite.

6.1. *An application.* We applied the above Exact Algorithms to the solution $X$ of the SDE (10) considered in Section 5. We compare the exact draws of our algorithms with the approximate ones of the simple Euler scheme that considers the continuous time process $Y = \{Y_t\}_{t \geq 0}$, $Y_0 = 0$, defined on the instances $\{ih\}_{i \geq 1}$, for some chosen increment $h > 0$, via the recursion

$$Y_{ih} = Y_{(i-1)h} + \sin(Y_{(i-1)h})h + Z_{i,h},$$

where $Z_{i,h}$, $i \geq 1$, are i.i.d. draws from the normal distribution with mean 0 and variance $h$. The paths of $Y$ become continuous after considering the linear interpolations between the successive instances of the grid $\{ih\}_{i \geq 1}$.

On the left of Figure 6, we show a *qq*-plot comparing two samples each of size 50,000 from the distribution of $M_2^X = \sup\{X_t; 0 \leq t \leq 2\}$ generated by the Exact Algorithm and the Euler approximation. For the Euler scheme we used increments $h = 2^{-9}$. On the right of the same figure, we show the times in seconds needed to get 50,000 draws from $M_2^X$ for the Exact Algorithm and the Euler scheme for different increments $h$. We have also included the results of a Kolmogorov–Smirnov test that compares the approximate samples of the Euler scheme with the exact sample as an indication of how small $h$ has to be for the Euler scheme to give correct results.

In Figure 7 we present similar results for the case of the simulation of the first passage time $\tau_2 = \inf\{t \geq 0 : X_t = 2\}$. Actually we tune the algorithm to obtain draws from the $\min\{\tau_2, 10\}$ since it can be shown that $\mathrm{E}[\tau_2] = \infty$ so the expected number of the uniformly drawn points for the completion of the algorithm for drawing from $\tau_2$ is infinite.

It is remarkable that in both cases the Exact Algorithm appears much more efficient than the Euler approximation.

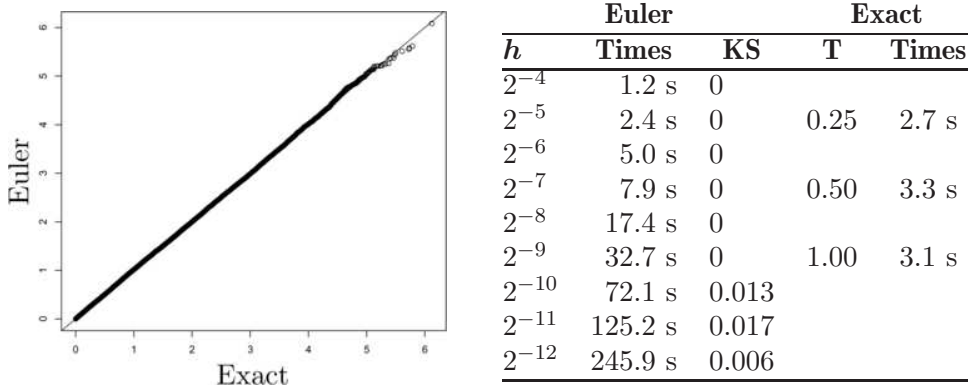

| | **Euler** | | **Exact** | |
|---|---:|---|---:|---:|
| $h$ | **Times** | **KS** | **T** | **Times** |
| $2^{-4}$ | 1.2 s | 0 | | |
| $2^{-5}$ | 2.4 s | 0 | 0.25 | 2.7 s |
| $2^{-6}$ | 5.0 s | 0 | | |
| $2^{-7}$ | 7.9 s | 0 | 0.50 | 3.3 s |
| $2^{-8}$ | 17.4 s | 0 | | |
| $2^{-9}$ | 32.7 s | 0 | 1.00 | 3.1 s |
| $2^{-10}$ | 72.1 s | 0.013 | | |
| $2^{-11}$ | 125.2 s | 0.017 | | |
| $2^{-12}$ | 245.9 s | 0.006 | | |

FIG. 6. *Results from the simulation of the maximum $M_2^X = \sup\{X_t; 0 \leq t \leq 2\}$ of the solution of* (10).

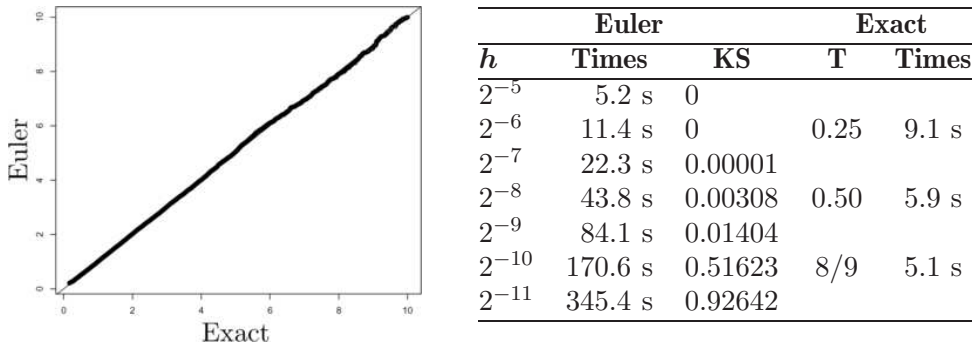

|       | Euler   |         |      | Exact  |
|-------|---------|---------|------|--------|
| $h$   | Times   | KS      | T    | Times  |
| $2^{-5}$  | 5.2 s   | 0       |      |        |
| $2^{-6}$  | 11.4 s  | 0       | 0.25 | 9.1 s  |
| $2^{-7}$  | 22.3 s  | 0.00001 |      |        |
| $2^{-8}$  | 43.8 s  | 0.00308 | 0.50 | 5.9 s  |
| $2^{-9}$  | 84.1 s  | 0.01404 |      |        |
| $2^{-10}$ | 170.6 s | 0.51623 | 8/9  | 5.1 s  |
| $2^{-11}$ | 345.4 s | 0.92642 |      |        |

FIG. 7. *Similar results to those of Figure 6 for the case of the simulation of $\min\{\tau_2, 10\}$ for $\tau_2 = \inf\{t \geq 0 : X_t = 2\}$.*

**7. Extensions and conclusions.** In this paper we have introduced a simple but computationally effective way of simulating exactly from a family of diffusion processes. The algorithm outputs a skeleton which can be readily "filled in" as and when necessary using simple Gaussian random variables, and crucially independently of the diffusion we are attempting to simulate from.

We have not carried out an extensive simulations study. However, in the examples considered, our method performs very favorably in comparison to the obvious numerical approximation alternative using the Euler scheme. In Section 4 we also give results which show that the method can be robust to the length of the time-series, and computing time is at worst linear in the degree of nonlinearity (as measured by the range of $\alpha^2 + \alpha'$).

The convenient form of the output allows the algorithm to be used in a number of ways including the construction of reduced variance Monte Carlo estimation. In Section 5 we discuss unbiased estimation of boundary hitting times and diffusion maxima. We envisage future application in inference for stochastic processes and finance.

The most demanding of the Conditions 1–3 (detailed in Section 3.1) required for the algorithm to work, is that the functional $\alpha^2 + \alpha'$ of the drift be bounded. In most cases $\alpha^2 + \alpha'$ is bounded from below but not from above so we can still produce (5) for some $\phi \geq 0$ and hope for a valid rejection sampling scheme. Ongoing work is investigating such an approach.

It is worth remarking that the approach outlined in this paper extends routinely to SDEs for non-Markov processes absolutely continuous with respect to appropriate martingales, and our focus on diffusions has been purely for simplicity. Furthermore, it is easy to extend these results some way toward considering jump diffusions.

APPENDIX

PROOF OF PROPOSITION 1 (Rejection sampling). Note that for any $i = 1, 2, \ldots$ we get

$$P[I_i = 1] = \int_S P[I_i = 1 | Y_i = y] \nu(dy) = \int_S f(y) \nu(dy) = \int_S \varepsilon \mu(dy) = \varepsilon.$$

Trivially, for any $F \in \mathcal{S}$ we have the property

(14) $\quad P[Y_\tau \in F] = P[Y_\tau \in F, I_1 = 1] + P[Y_\tau \in F | I_1 = 0] \cdot P[I_1 = 0].$

From the independence among the members of the sequence $(Y_n, I_n)_{n \geq 1}$ it is clear that

(15) $\quad\quad\quad\quad\quad P[Y_\tau \in F | I_1 = 0] = P[Y_\tau \in F].$

We can easily find the $P[Y_\tau \in F, I_1 = 1]$ in the following way:

(16) $\quad P[Y_\tau \in F, I_1 = 1] = \int_F P[I_1 = 1 | Y_1 = y] \nu(dy) = \int_F f(y) \nu(dy) = \varepsilon \mu(F).$

From (14) using (15) and (16) we get

$$P[Y_\tau \in F] = \varepsilon \mu(F) + (1 - \varepsilon) P[Y_\tau \in F],$$

which yields that $P[Y_\tau \in F] = \mu(F)$. □

PROOF OF PROPOSITION 3 (Biased Brownian motion). Choose any $F \in \mathcal{C}$. We will show that $E_\mathbb{Z}[\mathbb{I}_F] = E_\mathbb{W}[\mathbb{I}_F f]$ where we have set

$$f(\omega) := \frac{h(B_T)\sqrt{2\pi T}}{\exp(-B_T^2/(2T))}, \quad \omega \in C.$$

Note that $f$ is $\sigma(B_T)$-measurable. From the definition of $\overline{BM}$ it is clear that

$$\mathbb{Z}[F | \sigma(B_T)] = \mathbb{W}[F | \sigma(B_T)] =: g(B_T) \quad \mathbb{W}\text{-a.s.}$$

for some Borel-measurable function $g : \mathbf{R} \mapsto \mathbf{R}$. Clearly,

$$E_\mathbb{W}[\mathbb{I}_F f] = E_\mathbb{W}[E_\mathbb{W}[\mathbb{I}_F f | \sigma(B_T)]] = E_\mathbb{W}[f \mathbb{W}[F | \sigma(B_T)]]$$

$$= \int_\mathbf{R} \frac{h(u)\sqrt{2\pi T}}{\exp(-u^2/(2T))} \frac{\exp(-u^2/(2T))}{\sqrt{2\pi T}} g(u) \, du = \int_\mathbf{R} h(u) g(u) \, du$$

since w.r.t. $\mathbb{W}$ the random variable $B_T$ is distributed according to the normal distribution with mean 0 and variance $T$. It is straightforward that

$$E_\mathbb{Z}[\mathbb{I}_F] = E_\mathbb{Z}[E_\mathbb{Z}[\mathbb{I}_F | \sigma(B_T)]] = \int_\mathbf{R} h(u) g(u) \, du$$

since under $\mathbb{Z}$ we know that $B_T \sim h$. It is clear that $E_\mathbb{Z}[\mathbb{I}_F] = E_\mathbb{W}[\mathbb{I}_F f]$. □

PROOF OF THEOREM 1. (i) From (8) it is straightforward that $E_1 \subseteq E_3 \subseteq \cdots \subseteq E_{2n+1} \subseteq E_{2n+2}$ for any $n \in \{0,1,\ldots\}$. From the definition of $(E_{2n+2})_{n\geq 0}$ we can conclude that $E_{2n} = E_{2n+2} + (\Gamma_{2n+1} - \Gamma_{2n+2})$, so $E_{2n+2} \subseteq E_{2n} \subseteq \cdots \subseteq E_2$. Clearly, $E_{2\kappa+1} \subseteq E_{2\lambda+2}$ for any $\kappa, \lambda \in \{0,1,\ldots\}$.

The definition of $(\Gamma_n)_{n\geq 0}$ implies that for any $n \geq 0$ it is true that $\Gamma_{n+1} \subseteq \Gamma_n$, $\bigcap_n \Gamma_n$ is a set of zero probability and $\Gamma_{2n+2} = E_{2n+2} - E_{2n+1}$. Therefore:

$$\left(\bigcap_n E_{2n+2}\right) - \left(\bigcup_n E_{2n+1}\right) = \bigcap_n (E_{2n+2} - E_{2n+1}) = \bigcap_n \Gamma_{2n+2}.$$

Trivially, $\bigcup_n E_{2n+1} \subseteq \bigcap_n E_{2n+2}$ and their difference $\bigcap_n E_{2n+2} - \bigcup_n E_{2n+1}$ has zero probability.

(ii) Since $E_{2n+1} \uparrow E$ it is true that

$$\text{Prob}[I=1|\omega] \equiv \text{Prob}[E|\omega] = \lim_{n\to\infty} \text{Prob}[E_{2n+1}|\omega]. \tag{17}$$

Recall that $E_{2n+1} = \sum_{k=0}^n (\Gamma_{2k} - \Gamma_{2k+1})$. We can now get that

$$\begin{aligned}
\text{Prob}[E_{2n+1}|\omega] &= \sum_{k=0}^n (\text{Prob}[\Gamma_{2k}|\omega] - \text{Prob}[\Gamma_{2k+1}|\omega]) \\
&= \sum_{k=0}^{2n+1} (-1)^k \text{Prob}[\Gamma_k|\omega].
\end{aligned} \tag{18}$$

From the fact that $\tau = (V_n, W_n)_{n\geq 1}$ are i.i.d. and $\tau$, $U$ and $\omega$ are independent we conclude that

$$\text{Prob}[\Gamma_k|\omega] = \text{Prob}\left[U \leq \frac{1}{k!}\Big|\omega\right] \prod_{i=1}^k \text{Prob}[\phi(B_{V_i}(\omega)) \geq W_i|\omega]. \tag{19}$$

The $\text{Prob}[\phi(B_{V_i}(\omega)) \geq W_i|\omega]$ is the probability that given a path $\omega \sim \mathbb{Z}$ a point uniformly selected from the rectangle $(0,T) \times (0,1/T)$ is found below the graph $\{(t, \phi(B_t(\omega))) : t \in [0,T]\}$. Recall that $B_t$ is just the coordinate mapping $B_t = \omega(t)$, $t \in [0,T]$. From (6) and elementary probability theory we get

$$\text{Prob}[\phi(B_{V_k}(\omega)) \geq W_k|\omega] = \int_0^T \phi(B_t(\omega))\,dt.$$

It is clear that $\text{Prob}[U \leq \frac{1}{k!}|\omega] = \frac{1}{k!}$, so (19) yields

$$\text{Prob}[\Gamma_k|\omega] = \frac{1}{k!}\left\{\int_0^T \phi(B_t(\omega))\,dt\right\}^k. \tag{20}$$

Starting from (17) and using (18), (20) we get

$$\text{Prob}[I = 1|\omega] = \lim_{n \to \infty} \sum_{k=0}^{2n+1} \frac{(-1)^k}{k!} \left\{ \int_0^T \phi(B_t(\omega)) \, dt \right\}^k$$

$$= \exp\left\{ -\int_0^T \phi(B_t(\omega)) \, dt \right\}. \qquad \square$$

PROOF OF PROPOSITION 4. Since $0 \leq \phi \leq T^{-1}$, the probability of accepting an arbitrary path $\omega \sim \mathbb{Z}$ will be $\exp\{-\int_0^T \phi(B_t) \, dt\} \geq \exp(-1)$.

Recall that $T$ is restricted not to be bigger than some constant $T_0 (= \frac{1}{k_2 - k_1})$; see Condition 3 for the definition of $k_1, k_2$. To emphasize the involvement of the time variable in what follows we define $p : [0, T_0] \times C \mapsto [0, 1]$ with

$$p(T, \omega) = \exp\left\{ -\int_0^T \phi(B_t(\omega)) \, dt \right\}.$$

From result (ii) of Theorem 1 it is clear that

(21) $$\varepsilon(T) \equiv \int_C p(T, \omega) \mathbb{Z}(d\omega)$$

for all eligible $T$. Recall that $0 \leq \phi \leq T^{-1}$ for any $T \leq T_0$. Trivially, $p$ is decreasing in $T$ and $\lim_{T \downarrow 0} p(T, \omega) = 1$, both properties being true $\mathbb{Z}$-a.s. The first property yields, after using (21), that $\varepsilon(T)$ is decreasing in $T$. Also, since $\lim_{n \to \infty} p(1/n, \omega) = 1$ and $p(\frac{1}{n}, \omega) \leq p(\frac{1}{n+1}, \omega)$ for any $n \geq 1$ the monotone convergence theorem gives

$$\lim_{n \to \infty} \varepsilon(1/n) = \lim_{n \to \infty} \int_C p(1/n, \omega) \mathbb{Z}(d\omega) = \int_C \left\{ \lim_{n \to \infty} p(1/n, \omega) \right\} \mathbb{Z}(d\omega) = 1.$$

From the monotonicity of $\varepsilon(T)$ we conclude that $\lim_{T \downarrow 0} \varepsilon(T) = 1$. $\square$

PROOF OF PROPOSITION 6. For the time interval $[0, T]$, for $T = 1/R$, the algorithm uses $N_T := N_1 + N_2 + \cdots + N_J$ points, where $N_i$ is the number of points needed to decide about the $i$th proposed path and $J$ is the number of the proposed paths until one is accepted, $J = 1, 2, \ldots$ and $i = 1, 2, \ldots, J$. We denote by $N$ the expected number of points for deciding (in general) about a path $\omega \sim \mathbb{Z}$ and by $\mathrm{E}[N|A]$ and $\mathrm{E}[N|A^c]$ the expected number of points for deciding about a path $\omega \sim \mathbb{Z}$ given that the path is accepted and rejected, respectively. Let $\varepsilon = \text{Prob}[A]$ be the probability of accepting a path. Then

$$\mathrm{E}[N_T] = \mathrm{E}[\mathrm{E}[N_T|J]] = \mathrm{E}\left[ \sum_{i=1}^{J-1} \mathrm{E}[N_i|J] + \mathrm{E}[N_J|J] \right]$$

$$= \mathrm{E}[(J-1) \cdot \mathrm{E}[N|A^c] + \mathrm{E}[N|A]] = (1/\varepsilon - 1)\mathrm{E}[N|A^c] + \mathrm{E}[N|A]$$

$$= \frac{1}{\varepsilon} \cdot \mathrm{E}[N].$$

From Proposition 4 we have that $1/\varepsilon \leq e$ and from Proposition 5 that $E[N] \leq e$, so $E[N_T] \leq e^2$. The same result can be obtained for any $\lceil l \cdot R \rceil$ skeletons of length $T = 1/R$ (or less, for the case of the last skeleton when $1/R$ does not divide $l$) that will merge to construct the complete skeleton of length $l$. □

**Acknowledgments.** We thank Professor Duncan Murdoch for stimulating discussions on the problem of the exact simulation of diffusions. We also thank the referees for their useful suggestions.

## REFERENCES

[1] DEVROYE, L. (1986). *Nonuniform Random Variate Generation*. Springer, New York. MR836973
[2] FORSYTHE, G. E. (1972). von Neumann's comparison method for random sampling from the normal and other distributions. *Math. Comp.* **26** 817–826. MR315863
[3] KARATZAS, I. and SHREVE, S. E. (1991). *Brownian Motion and Stochastic Calculus*, 2nd ed. Springer, New York. MR1121940
[4] KENDALL, W. S. (1997). On some weighted Boolean models. In *Advances in Theory and Applications of Random Sets* (D. Jeulin, ed.) 105–120. World Scientific, Singapore. MR1654418
[5] KLOEDEN, P. E. and PLATEN, E. (1995). *Numerical Solution of Stochastic Differential Equations*. Springer, Berlin.
[6] LERCHE, H. R. (1986). *Boundary Crossing of Brownian Motion*. Springer, Berlin. MR861122
[7] MCKEAN, H. P. (1969). *Stochastic Integrals*. Academic Press, New York. MR247684
[8] OKSENDAL, B. K. (1998). *Stochastic Differential Equations*: *An Introduction with Applications*. Springer, Berlin. MR1619188
[9] PAPASPILIOPOULOS, O. and ROBERTS, G. O. (2005). Retrospective MCMC methods for Dirichlet process hierarchical models. Unpublished manuscript.

DEPARTMENT OF MATHEMATICS AND STATISTICS
LANCASTER UNIVERSITY
UNITED KINGDOM
E-MAIL: a.beskos@lancaster.ac.uk
        g.o.roberts@lancaster.ac.uk